\documentclass[12pt]{article}
\usepackage[reqno]{amsmath}
\usepackage{amssymb, theorem, enumerate, bm}
\usepackage{url}
\usepackage{hyperref}
\hypersetup{colorlinks=false}
\newcommand\cA{{\mathcal A}}

\newcommand\cC{{\mathcal C}}

\newcommand\cF{{\mathcal F}}
\newcommand\cG{{\mathcal G}}
\newcommand\cH{{\mathcal H}}

\newcommand\cS{{\mathcal S}}

\newcommand\cX{{\mathcal X}}
\newcommand\cY{{\mathcal Y}}

\makeatletter
\g@addto@macro\bfseries{\boldmath}
\makeatother

\newcommand\gbin[2]{\genfrac{[}{]}{0pt}{}{#1}{#2}}

\expandafter\ifx\csname pdfoptionalwaysusepdfpagebox\endcsname\relax\else
\fi

\theoremstyle{plain}
{\theorembodyfont{\slshape}
\newtheorem{theorem}{Theorem}[section]
\newtheorem{lemma}[theorem]{Lemma}
}
\theorembodyfont{\rmfamily}

\newtheorem{conjecture}[theorem]{Conjecture}

\newcommand\lref[1]{Lemma~\ref{lem:#1}}
\newcommand\tref[1]{Theorem~\ref{thm:#1}}
\newcommand\cref[1]{Corollary~\ref{cor:#1}}

\newcommand\conjref[1]{Conjecture~\ref{conj:#1}}

\def\sqr#1#2{{\vbox{\hrule height.#2pt
    \hbox{\vrule width.#2pt height#1pt \kern#1pt
        \vrule width.#2pt}\hrule height.#2pt}}}
\def\eqed{\sqr53}
\def\qed{%
    \ifmmode\eqno\eqed
    \else\nobreak\ \hfill\eqed\medbreak\fi}

\newcommand{\upperRomannumeral}[1]{\uppercase\expandafter{\romannumeral#1}}
\newcommand{\lowerromannumeral}[1]{\romannumeral#1\relax}

\title{The Manickam--Mikl\'{o}s--Singhi Conjectures for Sets and Vector Spaces}
\author{Ameera Chowdhury \thanks{Department of Mathematical Sciences, Carnegie Mellon University, Pittsburgh, PA, 15213, USA. E-mail: {\tt ameera@math.cmu.edu.} Research supported by NSF grant DMS-1203982.} \and Ghassan Sarkis \thanks{Department of Mathematics, Pomona College, Claremont, CA, 91711, USA. E-mail: {\tt Ghassan.Sarkis@pomona.edu, SShahriari@pomona.edu}} \and Shahriar Shahriari \footnotemark[2]}

\begin{document}
\maketitle

\begin{abstract}
More than twenty-five years ago, Manickam, Mikl\'{o}s, and Singhi conjectured that for positive integers $n,k$ with $n \geq 4k$, every set of $n$ real numbers with nonnegative sum has at least $\binom{n-1}{k-1}$ $k$-element subsets whose sum is also nonnegative. We verify this conjecture when $n \geq 8k^{2}$, which simultaneously improves and simplifies a bound of Alon, Huang, and Sudakov and also a bound of Pokrovskiy when $k < 10^{45}$.

Moreover, our arguments resolve the vector space analogue of this conjecture. Let $V$ be an $n$-dimensional vector space over a finite field. Assign a real-valued weight to each $1$-dimensional subspace in $V$ so that the sum of all weights is zero. Define the weight of a subspace $S \subset V$ to be the sum of the weights of all the $1$-dimensional subspaces it contains. We prove that if $n \geq 3k$, then the number of $k$-dimensional subspaces in $V$ with nonnegative weight is at least the number of $k$-dimensional subspaces in $V$ that contain a fixed $1$-dimensional subspace. This result verifies a conjecture of Manickam and Singhi from 1988.
\end{abstract}

\section{Introduction}
\label{sec:intro}

Manickam, Mikl\'{o}s, and Singhi \cite{mm, ms} conjectured in 1988 that
\begin{conjecture}
\label{conj:MMS}
For positive integers $n,k \in \mathbb{Z}^{+}$ with $n \geq 4k$, every set of $n$ real numbers with nonnegative sum has at least $\binom{n-1}{k-1}$ $k$-element subsets whose sum is also nonnegative.
\end{conjecture}
\conjref{MMS} was motivated by studies of the first distribution invariant in certain association schemes, and may also be considered an analogue of the Erd\H{o}s--Ko--Rado theorem \cite{ekr}. The Erd\H{o}s--Ko--Rado theorem states that if $n > 2k$, then any family of $k$-element subsets of an $n$-element set with the property that any two subsets have nonempty intersection has size at most $\binom{n-1}{k-1}$; moreover the unique extremal family is a star, the family of $k$-element subsets containing a fixed element.

\conjref{MMS} is similar to the Erd\H{o}s--Ko--Rado theorem, not only in the appearance of the binomial coefficient $\binom{n-1}{k-1}$, but also because the family of $k$-element subsets with nonnegative sum attains this lower bound and forms a star when one of the $n$ real numbers equals $n-1$ and the remaining $n-1$ numbers equal $-1$. As in the Erd\H{o}s--Ko--Rado theorem, $n$ must be large enough with respect to $k$, otherwise there exist $n$ real numbers with nonnegative sum and fewer than $\binom{n-1}{k-1}$ $k$-element subsets with nonnegative sum. Such examples can be easily constructed when $n=3k+r$ and $1 \leq r \leq k/7$. Although \conjref{MMS} and the Erd\H{o}s--Ko--Rado theorem share the same bound and extremal example, there is no obvious way to translate one question into the other.

\conjref{MMS} has attracted a lot of attention due to its connections with the Erd\H{o}s--Ko--Rado theorem \cite{AAH, ahs, aydblin, bhatt, bier, biermanickam, BC, blini, CIM, chowdhury, frankl, HS, HuangSudakov, manickam, mm, ms, MC, Pokrovskiy, Tyomkyn}, but still remains open. For more than two decades, \conjref{MMS} was known to hold only when $k|n$ \cite{ms} or when $n$ is at least an exponential function of $k$ \cite{bhatt, biermanickam, mm, Tyomkyn}. In their recent breakthrough paper, Alon, Huang, and Sudakov \cite{ahs} obtained the first polynomial bound $n \geq \min\{33k^{2}, 2k^{3}\}$ on \conjref{MMS}. Later, Aydinian and Blinovsky \cite{aydblin} and Frankl \cite{frankl} gave different proofs of \conjref{MMS} for a cubic range. Recently, a linear bound $n > 10^{46}k$ has been obtained by Pokrovskiy \cite{Pokrovskiy}. Finally, there are also several works that verify \conjref{MMS} for small $k$ \cite{chowdhury, HS, manickam, MC}.

\conjref{MMS} generalizes naturally to vector spaces. Let $V$ be an $n$-dimensional vector space over a finite field $\mathbb{F}_{q}$. For $k \in \mathbb{Z}^{+}$, we write $\gbin{V}{k}_{q}$ to denote the family of all $k$-dimensional subspaces of $V$. For $a, k \in \mathbb{Z}^{+}$, define the Gaussian binomial coefficient by
\begin{equation}
\label{Gaussiandef}
\gbin{a}{k}_{q}:= \prod_{0 \leq i < k} \frac{q^{a-i}-1}{q^{k-i}-1}.
\end{equation}
Simple counting arguments show that the number of $k$-dimensional subspaces in $V$ is $\gbin{n}{k}_{q}$ and that the number of $k$-dimensional subspaces in $V$ that contain a fixed $1$-dimensional subspace is $\gbin{n-1}{k-1}_{q}$. From now on, we omit the subscript $q$.

For each $1$-dimensional subspace $v \in \gbin{V}{1}$, assign a real-valued weight $f(v) \in \mathbb{R}$ so that the sum of all weights is zero. For a general subspace $S \subset V$, define its weight $f(S)$ to be the sum of the weights of all the $1$-dimensional subspaces it contains. We call a subspace \textsl{nonnegative} if it has nonnegative weight.

Manickam and Singhi posed the vector space analogue of \conjref{MMS} in 1988.
\begin{conjecture}[Manickam--Singhi, \cite{ms}]
\label{conj:vecanalog}
Let $V$ be an $n$-dimensional vector space over $\mathbb{F}_{q}$ and let $f: \gbin{V}{1} \rightarrow \mathbb{R}$ be a weighting of the 1-dimensional subspaces such that $\sum_{v \in \gbin{V}{1}} f(v) = 0$. If $n \geq 4k$, then there are at least $\gbin{n-1}{k-1}$ $k$-dimensional subspaces with nonnegative weight.
\end{conjecture}
Unlike \conjref{MMS}, we have no good reason for the $n \geq 4k$ stipulation as there are no known counterexamples to \conjref{vecanalog} for $n \geq 2k$. Counterexamples to \conjref{vecanalog} exist when $k < n < 2k$, and hence it is possible that \conjref{vecanalog} is true when $n \geq 2k$. Manickam and Singhi \cite{ms} showed that \conjref{vecanalog} holds when $k|n$.

\conjref{MMS} and \conjref{vecanalog} are strikingly similar, yet techniques that had been previously used to attack \conjref{MMS} do not readily generalize to vector spaces. As a simple example, note that the relative complement of a subset $S$ of $\{1, \ldots, n\}$ is both another subset of $\{1, \ldots, n\}$ and has empty intersection with $S$. In contrast, the relative complement of a subspace $A$ of $V$ is not another subspace of $V$, and the orthogonal complement of $A$ (with respect to an inner product) may have nontrivial intersection with $A$. Methods that work for vector spaces, however, do often
straightforwardly apply to sets and, for this reason, we are motivated to study \conjref{vecanalog}.

The main results of this paper verify \conjref{vecanalog} with a stronger statement and prove \conjref{MMS} when $n \geq 8k^{2}$. In particular, \tref{quadratic} simultaneously improves and simplifies the bound $n \geq \min\{33k^{2}, 2k^{3} \}$ of Alon, Huang, and Sudakov \cite{ahs} and also the bound $n \geq 10^{46}k$ of Pokrovskiy \cite{Pokrovskiy} when $k < 10^{45}$.

\begin{theorem}
\label{thm:vecanalog}
Let $V$ be an $n$-dimensional vector space over a finite field $\mathbb{F}_{q}$ and define $f: \gbin{V}{1} \rightarrow \mathbb{R}$ to be a weighting of the 1-dimensional subspaces such that $\sum_{v \in \gbin{V}{1}} f(v) = 0$. If $n \geq 3k$, then there are at least $\gbin{n-1}{k-1}$ $k$-dimensional subspaces with nonnegative weight. Moreover, if equality holds, the family of $k$-dimensional subspaces with nonnegative weight is a star, $\left \{S \in \gbin{V}{k} : \hat{v} \subset S \right \}$, on a fixed 1-dimensional subspace $\hat{v} \in \gbin{V}{1}$.
\end{theorem}

\begin{theorem}
\label{thm:quadratic}
Let $\cX = \{x_{1}, \ldots, x_{n} \} \subset \mathbb{R}$ be a set of $n$ real numbers whose sum is zero, and assume $x_{i} \geq x_{j}$ if $i \leq j$. If $n \geq 8k^{2}$, then at least $\binom{n-1}{k-1}$ $k$-element subsets of $\cX$ have nonnegative sum. Moreover, if equality holds, the family of $k$-element subsets with nonnegative sum is a star on $x_{1}$, $\left \{\cS \in \binom{\cX}{k} : x_{1} \in \cS \right \}$.
\end{theorem}
Note that there is no loss of generality in assuming that the $n$ real numbers in \tref{quadratic} are listed in decreasing order and sum to zero. Huang and Sudakov \cite{HuangSudakov} recently obtained similar results that are more general but which are weaker for both sets and vector spaces; in the latter case, Ihringer \cite{Ihringer} has recently extended our method to verify \conjref{vecanalog} for $n \geq 2k$ and large $q$.

\section{Bose-Mesner Matrices}
\label{sec:bosemesner}

The proofs of \tref{vecanalog} and \tref{quadratic} are similar. To avoid repetition and because the calculations in the vector space case are less familiar, we only give here the argument for the vector space case when the proof of the corresponding statement for sets is essentially the same. Full details for the case of sets are available in \cite{CSSQuad}.

We collect some notation and facts regarding the Gaussian binomial coefficients. When $k=1$, we write the Gaussian binomial coefficient $\gbin{a}{1}$ as $[a]$. A familiar relation involving binomial coefficients is Pascal's identity. We note two similar relations involving Gaussian binomial coefficients. For positive integers $a, k \in \mathbb{Z}^{+}$, we have
\begin{equation}
\label{qpascal}
\gbin{a}{k} = q^{a-k}\gbin{a-1}{k-1} + \gbin{a-1}{k} = \gbin{a-1}{k-1} +q^k \gbin{a-1}{k}.
\end{equation}

It will also be useful to note that given $S \in \gbin{V}{i}$ and $S' \in \gbin{V}{f}$ with $\dim(S \cap S') = 0$, the number of $e$-dimensional subspaces $T$ with $\dim(S' \cap T) = 0$ and $S \subset T$ is
\begin{equation}
\label{useful}
q^{f(e-i)} \gbin{n-i-f}{e-i}.
\end{equation}

We will need a lemma involving inclusion matrices $W_{jk}$ and Kneser matrices $\overline{W}_{jk}$. Let $W_{jk}$ (respectively $\overline{W}_{jk}$) denote the $\gbin{n}{j} \times \gbin{n}{k}$ matrix whose rows are indexed by the $j$-dimensional subspaces of $V$, whose columns are indexed by the $k$-dimensional subspaces of $V$, and where the entry in row $Y$ and column $S$ is $1$ if $Y \subset S$ (respectively if $Y \cap S = \langle 0 \rangle$) and is 0 otherwise. Similarly, let $\mathcal{W}_{jk}$ (respectively $\overline{\mathcal{W}}_{jk}$) denote the $\binom{n}{j} \times \binom{n}{k}$ matrix whose rows are indexed by the $j$-element subsets of $\cX$, whose columns are indexed by the $k$-element subsets of $\cX$, and where the entry in row $\cY$ and column $\cS$ is $1$ if $\cY \subset \cS$ (respectively if $\cY \cap \cS = \emptyset$) and is 0 otherwise.

Define $\vec{f} = (f(v_{1}), \ldots, f(v_{[n]})) \in \mathbb{R}^{[n]}$ to be a vector that gives the weights of each 1-dimensional subspace in \tref{vecanalog}. Similarly, define $\vec{x} = (x_{1}, \ldots, x_{n}) \in \mathbb{R}^{[n]}$ to be a vector that lists the $n$ real numbers in \tref{quadratic}. Without loss of generality, we may assume that $\vec{x}$ and $\vec{f}$ are nonzero, and observe that $\vec{x}$ and $\vec{f}$ are orthogonal to $\vec{1}$. The entries of $W_{1k}^{T} \vec{f}$ give the weights of $k$-dimensional subspaces of $V$. Similarly, the entries of $\mathcal{W}_{1k}^{T} \vec{x}$ give the sums of $k$-element subsets of $\cX$. We will show that $W_{1k}^{T} \vec{f}$ has at least $\gbin{n-1}{k-1}$ nonnegative entries when $n \geq 3k$. Similarly, we will show that $\mathcal{W}_{1k}^{T} \vec{x}$ has at least $\binom{n-1}{k-1}$ nonnegative entries when $n \geq 8k^{2}$.

An important observation by Frankl and Wilson \cite{FranklWilson} is that $W_{1k}^{T} \vec{f}$ is an eigenvector of the Bose-Mesner matrix
\begin{equation}
\label{bosemesner}
B_{j} = \overline{W}_{jk}^{T}W_{jk}
\end{equation}
for $0 \leq j \leq k$ with eigenvalue $-q^{j(k-1)} \gbin{k-1}{j-1} \gbin{n-j-1}{k-1}$. We include a proof for completeness.

\begin{lemma}[Frankl--Wilson, \cite{FranklWilson}]
\label{lem:eigenvec}
For $0 \leq j \leq k$, we have $W_{1k}^{T} \vec{f}$ is an eigenvector of the Bose-Mesner matrix $B_{j}$ with eigenvalue
\begin{equation}
\label{vectoreigenval}
-q^{j(k-1)} \gbin{k-1}{j-1} \gbin{n-j-1}{k-1}.
\end{equation}
\end{lemma}

\noindent \textbf{Proof.} Since the columns of $W_{1k}^{T}$ are linearly independent \cite{GJ, Kantor, WilsonNecessary} and $\vec{f} \neq \vec{0}$, we have that $W_{1k}^{T} \vec{f} \neq \vec{0}$. For $S \in \gbin{V}{j}$ and $T \in \gbin{V}{1}$, observe that the entry of $W_{jk}W_{1k}^{T}$ in row $S$ and column $T$ counts the number of $k$-dimensional subspaces of $V$ that contain $S$ and contain $T$. Hence,
\begin{equation}
\label{SunionT}
W_{jk}W_{1k}^{T}(S,T) = \begin{cases}
                        \gbin{n-j}{k-j} & \text{if $T \subset S$} \\
                        \gbin{n-j-1}{k-j-1} & \text{if $T \not \subset S$.} \\
                        \end{cases}
\end{equation}
For the remainder of this proof, $J$ is a matrix all of whose $\gbin{n}{1}$ columns are $\vec{1}$. Since $J \vec{f} = \vec{0}$, \eqref{qpascal} yields
\begin{equation}
\label{j1}
W_{jk}W_{1k}^{T}\vec{f} = \left( q^{k-j} \gbin{n-j-1}{k-j} W_{1j}^{T} + \gbin{n-j-1}{k-j-1}J \right)\vec{f} = q^{k-j} \gbin{n-j-1}{k-j}W_{1j}^{T} \vec{f}.
\end{equation}
For $A \in \gbin{V}{k}$ and $B \in \gbin{V}{1}$, observe that $\overline{W}_{jk}^{T} W_{1j}^{T}(A,B)$ counts the number of $j$-dimensional subspaces of $V$ that have trivial intersection with $A$ and that contain $B$. By \eqref{useful},
\begin{align}
\label{kk}
\overline{W}_{jk}^{T}W_{1j}^{T} \vec{f} &= q^{k(j-1)} \gbin{n-k-1}{j-1} \overline{W}_{1k}^{T} \vec{f} = q^{k(j-1)} \gbin{n-k-1}{j-1}(J - W_{1k}^{T}) \vec{f} \notag \\
&= -q^{k(j-1)} \gbin{n-k-1}{j-1} W_{1k}^{T} \vec{f},
\end{align}
since $J \vec{f} = \vec{0}$.
Multiplying \eqref{j1} on the left by $\overline{W}_{jk}^{T}$ and applying \eqref{kk} yields
\begin{align}
\label{eigenval}
B_{j}W_{1k}^{T} \vec{f} &= - q^{k(j-1)+(k-j)} \gbin{n-j-1}{k-j} \gbin{n-k-1}{j-1} W_{1k}^{T} \vec{f} \notag \\
&= -q^{j(k-1)} \gbin{k-1}{j-1} \gbin{n-j-1}{k-1} W_{1k}^{T} \vec{f},
\end{align}
which proves that $W_{1k}^{T} \vec{f}$ is an eigenvector of the Bose-Mesner matrix $B_{j}$ with eigenvalue $-q^{j(k-1)} \gbin{k-1}{j-1} \gbin{n-j-1}{k-1}$. \qed

A similar calculation shows that

\begin{lemma}[Wilson, \cite{WilsonExact}]
\label{lem:eigenvecset}
For $0 \leq j \leq k$, we have $\mathcal{W}_{1k}^{T} \vec{x}$ is an eigenvector of the Bose-Mesner matrix $\mathcal{B}_{j} = \overline{\mathcal{W}}_{jk}^{T} \mathcal{W}_{jk}$ with eigenvalue
\begin{equation}
\label{seteigenval}
-\binom{k-1}{j-1} \binom{n-j-1}{k-1}.
\end{equation}
\end{lemma}

\section{Bounds from Eigenvalues}
\label{sec:linear}

\lref{lotson1} and its set analogue \lref{lotson1sets} are the main results of this section. \lref{lotson1} shows that if $n \geq 3k+1$ and there are at most $\gbin{n-1}{k-1}$ nonnegative $k$-dimensional subspaces in $V$, then there exists a $1$-dimensional subspace $v \in \gbin{V}{1}$ with almost $\gbin{n-1}{k-1}$ nonnegative $k$-dimensional subspaces containing it. Similarly, \lref{lotson1sets} shows that if $n \geq Ck^{2}$ and there are at most $\binom{n-1}{k-1}$ nonnegative $k$-element subsets of $\cX$, then at least $(1 - \frac{6}{C}) \binom{n-1}{k-1}$ $k$-element subsets on $x_{1}$ have nonnegative sum.

Henceforth, we write $W_{1k}^{T} \vec{f} = \vec{b}$, and index the entries of $\vec{b}$ with $k$-dimensional subspaces $S \in \gbin{V}{k}$. Let $A \in \gbin{V}{k}$ denote the $k$-dimensional subspace of $V$ with highest weight. We first use \lref{eigenvec} to obtain a lower bound on the number of nonnegative $k$-dimensional subspaces that nontrivially intersect $A$.

\begin{lemma}
\label{lem:highestweightintersection}
There are greater than $q^{k(k-1)} \gbin{n-k-1}{k-1}$ nonnegative $k$-dimensional subspaces of $V$ that nontrivially intersect the highest weight $k$-dimensional subspace $A \in \gbin{V}{k}$.
\end{lemma}

\noindent \textbf{Proof.} Since $A$ is the highest weight $k$-dimensional subspace of $V$, we have $b_{A}$ is a largest entry of $\vec{b}$. Note that $b_{A} > 0$ since $\vec{b} \neq \vec{0}$ and $\vec{b}$ is orthogonal to $\vec{1}$.

For $S, T \in \gbin{V}{k}$, observe that $B_{j}(S,T)$ counts the number of $j$-dimensional subspaces of $V$ that lie in $T$ and have trivial intersection with $S$, so by \eqref{useful},
\begin{equation}
\label{bjentries}
B_{j}(S,T) = q^{j \dim(S \cap T)} \gbin{k- \dim(S \cap T)}{j}.
\end{equation}

By \lref{eigenvec} with $j=k$, the dot product of the row of $B_{k}$ corresponding to $A$ and $\vec{b}$ equals $-q^{k(k-1)} \gbin{n-k-1}{k-1} b_{A}$. Hence, \eqref{bjentries} yields
\begin{equation}
\label{a0}
\sum_{\dim(S \cap A) = 0} b_{S} = -q^{k(k-1)} \gbin{n-k-1}{k-1} b_{A}.
\end{equation}
Since $\vec{b}$ is orthogonal to $\vec{1}$, we see that
\begin{equation}
\label{flip}
\sum_{\dim(S \cap A) \neq 0} b_{S} = q^{k(k-1)} \gbin{n-k-1}{k-1} b_{A}.
\end{equation}
As $b_{A}$ is a largest entry of $\vec{b}$, there are greater than $q^{k(k-1)} \gbin{n-k-1}{k-1}$ $k$-dimensional subspaces of $V$ with nonnegative weight that nontrivially intersect $A$. \qed

Let $\mathcal{A} = \{x_{1}, \ldots, x_{k}\}$ and note that $\mathcal{A}$ is the $k$-element subset of $\cX$ with largest sum. Using \lref{eigenvecset}, we can obtain a set analogue of \lref{highestweightintersection} that gives a lower bound on the number of nonnegative $k$-element subsets of $\cX$ that intersect $\mathcal{A}$.

\begin{lemma}
\label{lem:highestweightintersectionset}
There are greater than $\binom{n-k-1}{k-1}$ $k$-element subsets of $\cX$ with nonnegative sum that intersect $\mathcal{A} = \{x_{1}, \ldots, x_{k} \}$.
\end{lemma}

We note two simple inequalities that will be useful in our computations.
\begin{align}
\frac{a-1}{b-1} &< \frac{a}{b}  &\mbox{for $b > a \geq 1$} \label{simple1} \\
\frac{[b]}{[a]} &< q^{b-a+1}    &\mbox{for $b \geq a \geq 1$ and $q \geq 2$.} \label{simple2}
\end{align}

We now prove a lemma that shows how many nonnegative $k$-dimensional subspaces \lref{highestweightintersection} guarantees with respect to $\gbin{n-1}{k-1}$.

\begin{lemma}
\label{lem:usefulcomputation}
For $k \leq a \leq n-k$, we have
\begin{equation}
\label{whatfraction}
q^{k(k-1)} \gbin{n-a-1}{k-1} \geq \frac{1}{q^{(a-k)(k-1)}} \left(1 - \frac{1}{q^{n-a-k}} \right) \gbin{n-1}{k-1}.
\end{equation}
\end{lemma}

\noindent \textbf{Proof.} By \eqref{useful}, we have that $q^{a(k-1)} \gbin{n-a-1}{k-1}$ counts the number of $k$-dimensional subspaces in $V$ that contain $S \in \gbin{V}{1}$ and have trivial intersection with $S' \in \gbin{V}{a}$ when $\dim(S \cap S') = 0$. Hence,
\begin{equation}
\label{disjoint}
q^{a(k-1)} \gbin{n-a-1}{k-1} \geq \gbin{n-1}{k-1} - [a] \gbin{n-2}{k-2},
\end{equation}
because $[a] \gbin{n-2}{k-2}$ is an upper bound on the number of $k$-dimensional subspaces in $V$ that contain $S \in \gbin{V}{1}$ and nontrivially intersect $S' \in \gbin{V}{a}$. Applying \eqref{simple1} and \eqref{simple2}, we have
\begin{equation}
\label{firstterm}
\frac{ [a] \gbin{n-2}{k-2} }{ \gbin{n-1}{k-1} } = \frac{[a][k-1]}{[n-1]} < \frac{1}{q^{n-a-k}}.
\end{equation}
Putting \eqref{disjoint} and \eqref{firstterm} together yields the lemma. \qed

Recall that $A \in \gbin{V}{k}$ is the highest weight $k$-dimensional subspace of $V$. Let $C \in \gbin{V}{k}$ denote the $k$-dimensional subspace of $V$ with highest weight such that $\dim(A \cap C) = 1$. We now use \lref{eigenvec} to obtain a lower bound on $b_{C}$, the weight of $C$, under the assumption that there are at most $\gbin{n-1}{k-1}$ $k$-dimensional subspaces with nonnegative weight in $V$.

\begin{lemma}
\label{lem:packing}
Let $A \in \gbin{V}{k}$ denote the highest weight $k$-dimensional subspace, and let $C \in \gbin{V}{k}$ denote the highest weight $k$-dimensional subspace such that $\dim(A \cap C) = 1$. If there are at most $\gbin{n-1}{k-1}$ nonnegative $k$-dimensional subspaces of $V$ then $b_{C}$, the weight of $C$, satisfies
\begin{equation}
\label{bc}
b_{C} \geq \left( 1 - \frac{q+1}{q^{n-2k+1}} \right) b_{A}.
\end{equation}
\end{lemma}

\noindent \textbf{Proof.} By \lref{eigenvec} with $j=k-1$, the dot product of the row of $B_{k-1}$ corresponding to $A$ and $\vec{b}$ equals $-q^{(k-1)^{2}}[k-1]\gbin{n-k}{k-1}b_{A}$. Hence, by  \eqref{bjentries},
\begin{equation}
\label{bk-1}
[k] \sum_{\dim(S \cap A) = 0} b_{S} + q^{k-1} \sum_{\dim(S \cap A) = 1} b_{S} = -q^{(k-1)^{2}}[k-1]\gbin{n-k}{k-1}b_{A}.
\end{equation}
By \eqref{qpascal}, \eqref{a0}, and \eqref{bk-1},
\begin{align}
\label{a1}
\sum_{\dim(S \cap A) = 1} b_{S} &= \frac{1}{q^{k-1}} \left(q^{k(k-1)} [k] \gbin{n-k-1}{k-1} - q^{(k-1)^{2}}[k-1]\gbin{n-k}{k-1} \right) b_{A} \notag \\
&= \left( q^{k(k-1)} \gbin{n-k-1}{k-1} - q^{(k-1)(k-2)}[k-1]\gbin{n-k-1}{k-2} \right)b_{A}.
\end{align}

By definition, $C \in \gbin{V}{k}$ is the $k$-dimensional subspace with highest weight that satisfies $\dim(A \cap C) = 1$. We claim that if there are at most $\gbin{n-1}{k-1}$ nonnegative $k$-dimensional subspaces then
\begin{equation}
\label{packingweight}
b_{C} \geq \frac{q^{k(k-1)} \gbin{n-k-1}{k-1} - q^{(k-1)(k-2)}[k-1]\gbin{n-k-1}{k-2}}{\gbin{n-1}{k-1}} b_{A}.
\end{equation}
We must have that \eqref{packingweight} holds, otherwise as $b_{A}$ is a largest entry of $\vec{b}$, \eqref{a1} implies there are greater than $\gbin{n-1}{k-1}$ nonnegative entries $b_{S}$ where $\dim(S \cap A) = 1$.

Now, we show that the fraction on the right hand of \eqref{packingweight} is at least the fraction on the right hand side of \eqref{bc}. Applying \eqref{simple1} and \eqref{simple2}, we have
\begin{align}
\label{secondterm}
\frac{ q^{(k-1)(k-2)} [k-1] \gbin{n-k-1}{k-2} }{ \gbin{n-1}{k-1} } & < \frac{ q^{(k-1)^{2}} \gbin{n-k-1}{k-2} }{ \gbin{n-1}{k-1} } \notag \\
&= q^{(k-1)^{2}} \frac{ [k-1] [n-k] \cdots [n-2k+2] } {[n-1][n-2] \cdots [n-k]} \notag \\
&< \frac{ q^{(k-1)^{2}} }{ q^{n-k} q^{(k-2)(k-1)}}  = \frac{1}{q^{n-2k+1}}.
\end{align}
Applying \lref{usefulcomputation} and \eqref{secondterm} yields the lemma. \qed

Recall that $\mathcal{A} = \{x_{1}, \ldots, x_{k} \}$ is the $k$-element subset of $\cX$ with largest sum. Let $\mathcal{C} = \{x_{1}, x_{k+1}, \ldots, x_{2k-1} \}$ and note that $\mathcal{C}$ is the $k$-element subset of $\cX$ with largest sum such that $|\mathcal{A} \cap \mathcal{C}| = 1$. Henceforth, we write $\mathcal{W}_{1k}^{T} \vec{x} = \vec{\mathit{b}}$, and index the entries of $\vec{\mathit{b}}$ with $k$-element subsets $\cS \in \binom{\cX}{k}$. Using \lref{eigenvecset} we can obtain a set analogue of \lref{packing} that gives a lower bound on $\mathit{b}_{\mathcal{C}}$, the sum of $\mathcal{C}$, under the assumptions that $n \geq k^{2}$ and that there are at most $\binom{n-1}{k-1}$ $k$-element subsets with nonnegative sum in $\cX$.

\begin{lemma}
\label{lem:packingset}
Let $\mathcal{A} = \{x_{1}, \ldots, x_{k} \}$ and let $\mathcal{C} = \{x_{1}, x_{k+1}, \ldots, x_{2k-1} \}$. If $n \geq k^{2}$ and there are at most $\binom{n-1}{k-1}$ nonnegative $k$-element subsets of $\cX$ then $\mathit{b}_{\mathcal{C}}$, the sum of $\mathcal{C}$, satisfies
\begin{equation}
\label{bcset}
\mathit{b}_{\mathcal{C}} \geq \left( 1 - \frac{(2k-1)(k-1)}{n-2k+1} \right) \mathit{b}_{\mathcal{A}}.
\end{equation}
\end{lemma}

\noindent \textbf{Proof.} Calculations similar to those in the proof of \lref{packing} show that
\begin{equation}
\label{a1set}
\sum_{|\cS \cap \cA| = 1} b_{\cS} = \left( \binom{n-k-1}{k-1} -(k-1)\binom{n-k-1}{k-2} \right) b_{\cA},
\end{equation}
which is nonnegative exactly when $n \geq k^{2}$. As a result, arguments similar to those in \lref{packing} yield that
\begin{equation}
\label{similartopacking}
\mathit{b}_{\mathcal{C}} \geq \left( \frac{ \binom{n-k-1}{k-1} - (k-1)\binom{n-k-1}{k-2} }{ \binom{n-1}{k-1} } \right) \mathit{b}_{\mathcal{A}}.
\end{equation}
Now, we show that the fraction on the right hand of \eqref{similartopacking} is at least the fraction on the right hand side of \eqref{bcset}. We have
\begin{equation}
\label{bcsetsimplify1}
\binom{n-k-1}{k-1} - (k-1) \binom{n-k-1}{k-2} = \left( 1 - \frac{(k-1)^{2}}{n-2k+1} \right) \binom{n-k-1}{k-1}.
\end{equation}
We also have that
\begin{equation}
\label{bcsetsimplify2}
\frac{ \binom{n-k-1}{k-1} }{ \binom{n-1}{k-1} } = \frac{ (n-k-1) \cdots (n-2k+1) }{ (n-1) \cdots (n-k+1) } > \left( 1 - \frac{k}{n-k+1} \right)^{k-1} > 1 - \frac{k(k-1)}{n-k+1}.
\end{equation}
Putting \eqref{bcsetsimplify1} and \eqref{bcsetsimplify2} together yields \eqref{bcset}. \qed

Suppose that $A \cap C = v \in \gbin{V}{1}$. We count the number of $k$-dimensional subspaces of $V$ that nontrivially intersect both $A$ and $C$ but do not contain $v$.

\begin{lemma}
\label{lem:small}
Let $A \in \gbin{V}{k}$ denote the highest weight $k$-dimensional subspace, and let $C \in \gbin{V}{k}$ denote the highest weight $k$-dimensional subspace such that $\dim(A \cap C) = 1$. Suppose that $A \cap C = v \in \gbin{V}{1}$. The number of $k$-dimensional subspaces of $V$ that nontrivially intersect both $A$ and $C$ but do not contain $v$ is at most
\begin{equation}
\label{constant}
\frac{1}{q^{n-3k}} \gbin{n-1}{k-1}.
\end{equation}
In particular, when $k=2$, the number of $2$-dimensional subspaces of $V$ that nontrivially intersect both $A$ and $C$ but do not contain $v$ is $q^{2}$.
\end{lemma}

\noindent \textbf{Proof.} A $k$-dimensional subspace that nontrivially intersects $A$ and $C$ but does not contain $v$ must contain a $2$-dimensional subspace in $A \vee C$ that intersects each of $A$ and $C$ in exactly one 1-dimensional subspace not equal to $v$. We now show that the number of such $2$-dimensional subspaces is
\begin{equation}
\label{additionalimprove}
([k]-1)([2k-2] - q^{k}[k-2] - 1) = ([k]-1)^{2}.
\end{equation}
For a $1$-dimensional subspace $w \in A$ not equal to $v$, the number of $2$-dimensional subspaces in $A \vee C$ that contain $w$ is $[2k-2]$ because $\dim(A \vee C) = 2k-1$. Of these $2$-dimensional subspaces containing $w$, we have that $q^{k}[k-2]+1$ of them lie in $A$ or have trivial intersection with $C$ by \eqref{useful}. For each of the 2-dimensional subspaces counted by \eqref{additionalimprove}, there are $\gbin{n-2}{k-2}$ $k$-dimensional subspaces which contain that $2$-dimensional subspace. Applying \eqref{simple1} and \eqref{simple2}, the number of $k$-dimensional subspaces that nontrivially intersect $A$ and $C$ but do not contain $v$ is at most
\begin{equation}
\label{ahhh}
([k]-1)^{2} \gbin{n-2}{k-2} < q^{2k} \frac{ [k-1] }{ [n-1] } \gbin{n-1}{k-1} < \frac{1}{q^{n-3k}} \gbin{n-1}{k-1}.
\end{equation}
In particular, when $k=2$, we have that \eqref{additionalimprove} equals $q^{2}$. \qed

Recall that $A \cap C = v \in \gbin{V}{1}$. Now we use \lref{eigenvec} to obtain a lower bound on the number of nonnegative $k$-dimensional subspaces that contain $v$ under the assumption that there are at most $\gbin{n-1}{k-1}$ $k$-dimensional subspaces of $V$ with nonnegative weight.

\begin{lemma}
\label{lem:lotson1}
Let $A \in \gbin{V}{k}$ denote the highest weight $k$-dimensional subspace, and let $C \in \gbin{V}{k}$ denote the highest weight $k$-dimensional subspace such that $\dim(A \cap C) = 1$. Suppose that $A \cap C = v \in \gbin{V}{1}$. If there are at most $\gbin{n-1}{k-1}$ nonnegative $k$-dimensional subspaces of $V$, then the number of nonnegative $k$-dimensional subspaces that contain $v$ is at least
\begin{equation}
\label{generallotson1}
\left( 1 - \frac{1}{q^{n-3k}} - \frac{1}{q^{n-2k-1}} - \frac{1}{q^{n-2k}} - \frac{1}{q^{n-2k+1}} + \frac{q+1}{q^{2n-4k+1}} \right) \gbin{n-1}{k-1}.
\end{equation}
In particular, when $k=2$, the number of nonnegative $2$-dimensional subspaces that contain $v$ is at least
\begin{equation}
\label{lotson1k2}
\left( 1 - \frac{1}{q^{n-6}} - \frac{1}{q^{n-3}} + \frac{q+1}{q^{2n-7}} \right) [n-1].
\end{equation}
\end{lemma}

\noindent \textbf{Proof.} By \lref{eigenvec} with $j=k$, the dot product of the row of $B_{k}$ corresponding to $C$ and $\vec{b}$ equals $-q^{k(k-1)}\gbin{n-k-1}{k-1} b_{C}$. Hence, \eqref{bjentries} with $j=k$ and \lref{packing} yield that
\begin{align}
\label{rinseandrepeat}
\sum_{\substack{ \dim(S \cap C) \neq 0, \\ \dim(S \cap A) \neq 0}} b_{S} + \sum_{\substack{\dim(S \cap C) \neq 0, \\ \dim(S \cap A) = 0}} b_{S} &= \sum_{\dim(S \cap C) \neq 0} b_{S} = q^{k(k-1)}\gbin{n-k-1}{k-1} b_{C} \notag \\
&\geq q^{k(k-1)} \left(1 - \frac{q+1}{q^{n-2k+1}} \right) \gbin{n-k-1}{k-1} b_{A}.
\end{align}
We claim that
\begin{equation}
\label{nottoobig}
\sum_{\substack{\dim(S \cap C) \neq 0, \\ \dim(S \cap A) = 0}} b_{S} \leq \left( \gbin{n-1}{k-1} - q^{k(k-1)} \gbin{n-k-1}{k-1} \right) b_{A}.
\end{equation}
If \eqref{nottoobig} does not hold, there would be at least $\gbin{n-1}{k-1} - q^{k(k-1)} \gbin{n-k-1}{k-1}$ nonnegative entries $b_{S}$ such that $\dim(S \cap C) \neq 0$ and $\dim(S \cap A) = 0$. By \lref{highestweightintersection}, this would yield at least $\gbin{n-1}{k-1} + 1$ nonnegative $k$-dimensional subspaces of $V$. Hence, \eqref{nottoobig} holds. Since $q \geq 2$, \lref{usefulcomputation} implies that
\begin{align}
\label{prettybig}
\sum_{\substack{\dim(S \cap C) \neq 0, \\ \dim(S \cap A) \neq 0}} b_{S} &\geq \left( q^{k(k-1)} \left(2 - \frac{q+1}{q^{n-2k+1}} \right) \gbin{n-k-1}{k-1} - \gbin{n-1}{k-1} \right) b_{A} \notag \\
&\geq \left( \left(2 - \frac{q+1}{q^{n-2k+1}} \right) \left( 1 - \frac{1}{q^{n-2k}} \right) - 1 \right) \gbin{n-1}{k-1} b_{A} \notag \\
%&= \left( 1 - \frac{q+1}{q^{n-2k+1}} - \frac{2}{q^{n-2k}} + \frac{q+1}{q^{2n-4k+1}} \right) \gbin{n-1}{k-1} b_{A} \notag \\
&\geq \left( 1 - \frac{1}{q^{n-2k-1}} - \frac{1}{q^{n-2k}} - \frac{1}{q^{n-2k+1}} + \frac{q+1}{q^{2n-4k+1}} \right) \gbin{n-1}{k-1} b_{A}.
\end{align}

By \lref{small}, at most $q^{-(n-3k)} \gbin{n-1}{k-1}$ subspaces $S \in \gbin{V}{k}$ nontrivially intersect $A$ and $C$ and do not contain $v$, and each such subspace has weight at most $b_{A}$. Hence, by \eqref{constant} and \eqref{prettybig},
\begin{align}
\label{improvonv}
&\sum_{v \subset S} b_{S} = \sum_{\substack{\dim(S \cap C) \neq 0, \\ \dim(S \cap A) \neq 0}} b_{S} - \sum_{\substack{\dim(S \cap C) \neq 0, \\ \dim(S \cap A) \neq 0, \\ v \not \subset S}} b_{S} \notag \\
%&\geq \left( 1 - \frac{1}{q^{n-2k-1}} - \frac{1}{q^{n-2k}} - \frac{1}{q^{n-2k+1}} + \frac{q+1}{q^{2n-4k+1}} \right)  \gbin{n-1}{k-1} b_{A} - \frac{1}{q^{n-3k}} \gbin{n-1}{k-1} b_{A} \notag \\
&\geq \left( 1 - \frac{1}{q^{n-3k}} - \frac{1}{q^{n-2k-1}} - \frac{1}{q^{n-2k}} - \frac{1}{q^{n-2k+1}}  + \frac{q+1}{q^{2n-4k+1}} \right) \gbin{n-1}{k-1} b_{A}.
\end{align}
When $k=2$, at most $q^{2} < q^{-(n-4)}[n-1]$ subspaces $S \in \gbin{V}{2}$ nontrivially intersect $A$ and $C$ and do not contain $v$ so \eqref{improvonv} can be improved. Hence, a lower bound on the number of nonnegative $k$-dimensional subspaces that contain $v$ is given by \eqref{generallotson1} and \eqref{lotson1k2} for $k \geq 3$ and $k=2$ respectively. \qed

Observe that in \eqref{improvonv}, we use the bound
\begin{equation}
\label{crudevec}
\sum_{\substack{\dim(S \cap C) \neq 0, \\ \dim(S \cap A) \neq 0, \\ v \not \subset S}} b_{S} \leq \frac{1}{q^{n-3k}} \gbin{n-1}{k-1} b_{A},
\end{equation}
which is not optimal because we bound by $b_{A}$ the weights of all $k$-dimensional subspaces $S \in \gbin{V}{k}$ which nontrivially intersect $A$ and $C$ and do not contain $v$. In the case of sets, we can apply more sophisticated counting to yield a better set analogue of \eqref{crudevec}. Now we prove a set analogue of \lref{lotson1} that gives a lower bound on the number of nonnegative $k$-element subsets that contain $x_{1}$ under the assumptions that $n \geq k^{2}$ and that there are at most $\binom{n-1}{k-1}$ $k$-element subsets of $\cX$ with nonnegative sum.

\begin{lemma}
\label{lem:lotson1sets}
If $n \geq k^{2}$ and there are at most $\binom{n-1}{k-1}$ nonnegative $k$-element subsets of $\cX$, then the number of nonnegative $k$-element subsets that contain $x_{1}$ is at least
\begin{equation}
\label{generallotson1set}
\left( 1 - \frac{(6k-3)(k-1)}{n-2k+1} \right) \binom{n-1}{k-1}.
\end{equation}
\end{lemma}

\noindent \textbf{Proof.} Recall that $\mathcal{A} = \{x_{1}, \ldots, x_{k} \}$ and that $\mathcal{C} = \{x_{1}, x_{k+1}, \ldots, x_{2k-1} \}$. \lref{packingset}, \eqref{bcsetsimplify2}, and calculations similar to those in \lref{lotson1} show that
\begin{equation}
\label{prettybigset}
\sum_{\substack{\cS \cap \cC \neq \emptyset, \\ \cS \cap \cA \neq \emptyset}} \mathit{b}_{\cS} \geq \left( 1 - \frac{(4k-1)(k-1)}{n-2k+1} \right) \binom{n-1}{k-1} \mathit{b}_{\cA}.
\end{equation}

Let $\bm{\cF}_{i}$ be the family of $k$-element subsets of $\cX$ that contain $x_{i}$ but not $x_{1}$ and intersect $\cA$ and $\cC$. We have
\begin{equation}
\label{whereweregoing}
\sum_{x_{1} \in \cS} \mathit{b}_{\cS} = \sum_{\substack{\cS \cap \cC \neq \emptyset, \\ \cS \cap \cA \neq \emptyset}} \mathit{b}_{\cS} - \sum_{\substack{\cS \cap \cC \neq \emptyset, \\ \cS \cap \cA \neq \emptyset, \\ x_{1} \notin \cS}} \mathit{b}_{\cS} = \sum_{\substack{\cS \cap \cC \neq \emptyset, \\ \cS \cap \cA \neq \emptyset}} \mathit{b}_{\cS} - \sum_{i=2}^{n} |\bm{\cF}_{i}|x_{i}.
\end{equation}
We first show that if $i \in \{2, \ldots, 2k-1\}$ then
\begin{equation}
\label{case1}
|\bm{\cF}_{i}| = \binom{n-2}{k-1} - \binom{n-k-1}{k-1}.
\end{equation}
Without loss of generality suppose that $x_{i} \in \cA \setminus \{x_{1}\}$. There are $\binom{n-2}{k-1}$ $k$-element sets of $\cX$ that contain $x_{i}$ but not $x_{1}$. From these, we subtract the $\binom{n-k-1}{k-1}$ $k$-element subsets of $\cX$ that contain $x_{i}$ but do not intersect $\cC$.

Now we determine $|\bm{\cF}_{i}|$ when $i \in \{2k, \ldots, n\}$. Let $\bm{\mathcal{G}}_{i}$ (respectively $\bm{\cH}_{i}$) be the family of $k$-element subsets of $\cX$ that contain $x_{i}$ but not $x_{1}$ and intersect $\cA$ (respectively $\cC$). We have $\bm{\cF}_{i} = \bm{\cG}_{i} \cap \bm{\cH}_{i}$ so by inclusion-exclusion,
\begin{align}
\label{inclexcl}
| \bm{\cF}_{i} | &= |\bm{\cG}_{i} \cap \bm{\cH}_{i}| = |\bm{\cG}_{i}| + |\bm{\cH}_{i}| - |\bm{\cG}_{i} \cup \bm{\cH}_{i}| \notag \\
&= 2 \left( \binom{n-2}{k-1} - \binom{n-k-1}{k-1} \right) - \left( \binom{n-2}{k-1} - \binom{n-2k}{k-1} \right) \notag \\
&= \binom{n-2}{k-1} - 2 \binom{n-k-1}{k-1} + \binom{n-2k}{k-1}.
\end{align}

By \eqref{whereweregoing}, \eqref{case1}, and \eqref{inclexcl},
\begin{align}
\label{evenmorebetter}
\sum_{\substack{\cS \cap \cC \neq \emptyset, \\ \cS \cap \cA \neq \emptyset, \\ x_{1} \notin \cS}} \mathit{b}_{\cS} &= |\bm{\cF}_{2}| \sum_{i=2}^{2k-1} x_{i} + |\bm{\cF}_{2k}| \sum_{i=2k}^{n} x_{i} = |\bm{\cF}_{2}|(\mathit{b}_{\cA} + \mathit{b}_{\cC} - 2x_{1}) + |\bm{\cF}_{2k}|(x_{1} - \mathit{b}_{\cA} - \mathit{b}_{\cC}) \notag \\
&= (2|\bm{\cF}_{2}| - |\bm{\cF}_{2k}|)(-x_{1}) + (|\bm{\cF}_{2}| - |\bm{\cF}_{2k}|)(\mathit{b}_{\cA} + \mathit{b}_{\cC}) \notag \\
&< 2(|\bm{\cF}_{2}| - |\bm{\cF}_{2k}|) \mathit{b}_{\cA} = 2 \sum_{j=k+2}^{2k} \binom{n-j}{k-2} \mathit{b}_{\cA} \notag \\
&< 2(k-1) \binom{n-k-2}{k-2} \mathit{b}_{\cA} < \frac{2(k-1)^{2}}{n-1} \binom{n-1}{k-1} \mathit{b}_{\cA}.
\end{align}

By \eqref{prettybigset}, \eqref{whereweregoing}, and \eqref{evenmorebetter}, we have
\begin{equation}
\label{almostthereset}
\sum_{x_{1} \in \cS} \mathit{b}_{\cS} \geq \left( 1 - \frac{(6k-3)(k-1)}{n-2k+1} \right) \binom{n-1}{k-1} b_{\cA}.
\end{equation}
Hence, a lower bound on the number of nonnegative $k$-element subsets that contain $x_{1}$ is given by \eqref{generallotson1set}. \qed

\section{Bounds from Averaging}
\label{sec:greenekleitman}

\lref{nonnegT} and its set analogue \lref{nonnegTset} are the main results of this section. \lref{nonnegT} shows that if $n \geq 3k$ and $T \in \gbin{V}{k}$ is a $k$-dimensional subspace of $V$ with negative weight, then there are almost $\gbin{n-1}{k-1}$ nonnegative $k$-dimensional subspaces that have trivial intersection with $T$. Similarly, \lref{nonnegTset} shows that if $\mathcal{T} \in \binom{\cX}{k}$ is a $k$-element subset with negative sum, then there are at least $\binom{n-2k}{k-1}$ nonnegative $k$-element subsets of $\cX$ that are disjoint from $\mathcal{T}$.

We first show that if $T \in \gbin{V}{k}$ has negative weight, then there is a $(k+i)$-dimensional subspace $W \in \gbin{V}{k+i}$ with negative weight that contains $T$ when $0 \leq i \leq n-k-1$.

\begin{lemma}
\label{lem:induct}
If $T \in \gbin{V}{k}$ has negative weight and $0 \leq i \leq n-k-1$, then there is a negative weight $(k+i)$-dimensional subspace $W \in \gbin{V}{k+i}$ containing $T$.
\end{lemma}

\noindent \textbf{Proof.} We induct on $0 \leq j < n-k-1$ to show that if there is a $(k+j)$-dimensional subspace $W_{j}$ with negative weight that contains $T$, then there is a $(k+j+1)$-dimensional subspace with negative weight that contains $T$. The hypothesis is true for $j=0$ since $T \in \gbin{V}{k}$ has negative weight, so we can set $W_{0} = T$. Let $W_{j+1,1}, \ldots, W_{j+1,[n-k-j]}$ denote all the $(k+j+1)$-dimensional subspaces containing $W_{j}$. Since $W_{j}$ has negative weight,
\begin{align}
\label{negtoneg}
\sum_{l=1}^{[n-k-j]} f(W_{j+1,l}) &= [n-k-j]f(W_{j}) + \sum_{w \in \gbin{V}{1} \setminus \gbin{W_{j}}{1}} f(w) = [n-k-j]f(W_{j}) - f(W_{j}) \notag \\
&= ([n-k-j]-1)f(W_{j}) < 0,
\end{align}
as $0 \leq j < n-k-1$. Hence, at least one of the $(k+j+1)$-dimensional subspaces $W_{j+1,l}$ counted in \eqref{negtoneg} must have negative weight. \qed

To find nonnegative $k$-dimensional subspaces that have trivial intersection with $T$, we introduce the concept of a $k$-spread. A family $\bm{S} \subset \gbin{V}{k}$ of $k$-dimensional subspaces is called a $k$-\textsl{spread} if every $1$-dimensional subspace of $V$ is contained in exactly one $k$-dimensional subspace in $\bm{S}$. Andr\'{e} \cite{andre} showed that a necessary and sufficient condition for a $k$-spread to exist in $V$ is that $k|n$. We need the following lemmas due to Beutelspacher \cite{mrB} that show how to construct a large partial $k$-spread in the case that $n$ is not divisible by $k$.

\begin{lemma}[Beutelspacher, \cite{mrB}]
\label{lem:mrB1}
Let $W$ be an $(s+t)$-dimensional vector space over $\mathbb{F}_{q}$, where $s \geq t$. If $U \in \gbin{W}{s}$ is an $s$-dimensional subspace of $W$, then there exists a family $\bm{S} \subset \gbin{W}{t}$ of $t$-dimensional subspaces of $W$ which have trivial intersection with $U$ and such that every $1$-dimensional subspace in $\gbin{W}{1} \setminus \gbin{U}{1}$ is contained in exactly one element of $\bm{S}$.
\end{lemma}

\noindent \textbf{Proof.} Let $Y$ be a $2s$-dimensional vector space over $\mathbb{F}_{q}$ that contains $W$, and let $\bm{S'}$ be an $s$-spread of $Y$ which contains $U$. Since each element of $\bm{S'} \setminus \{U\}$ has trivial intersection with $U$, we see that each element of $\bm{S'} \setminus \{U\}$ intersects $W$ in a $t$-dimensional subspace. The family of intersections $\bm{S} = \{ S \cap W : S \in \bm{S'} \setminus \{U\} \}$ has the required properties. \qed

Now we use \lref{mrB1} to construct large partial $k$-spreads in the case that $n$ is not divisible by $k$.

\begin{lemma}[Beutelspacher, \cite{mrB}]
\label{lem:mrBtotherescue}
Write $n = mk+r$ where $0 \leq r \leq k-1$. If $U \in \gbin{V}{k+r}$ is a $(k+r)$-dimensional subspace of $V$, then there exists a family $\bm{S} \subset \gbin{V}{k}$ of $k$-dimensional subspaces which
\begin{enumerate}[(i)]
\item have trivial intersection with $U$
\item and such that every $1$-dimensional subspace in $\gbin{V}{1} \setminus \gbin{U}{1}$ is contained in exactly one element of $\bm{S}$.
\end{enumerate}
Moreover,
\begin{equation}
\label{sizeS}
|\bm{S}| = \frac{q^{k+r}[n-k-r]}{[k]}.
\end{equation}
\end{lemma}

\noindent \textbf{Proof.} Apply \lref{mrB1} with $t=k$ and $s = k+r, 2k+r, \ldots, (m-1)k+r$ to construct $\bm{S}$. Since $\bm{S} \subset \gbin{V}{k}$ satisfies properties (\lowerromannumeral{1}) and (\lowerromannumeral{2}), the size of $\bm{S}$ is
\begin{equation}
\label{sizeSproof}
|\bm{S}| = \frac{[n]-[k+r]}{[k]} = \frac{ q^{n} - q^{k+r} }{q^{k} - 1} = \frac{q^{k+r}(q^{n-k-r}-1)}{q^{k}-1} = \frac{q^{k+r}[n-k-r]}{[k]}. \; \; \eqed
\end{equation}

Let $n = mk+r$, where $0 \leq r \leq k-1$, and suppose that $U \in \gbin{V}{k+r}$ is a negative weight $(k+r)$-dimensional subspace. \lref{permute} provides a way of finding $k$-dimensional subspaces with nonnegative weight that have trivial intersection with $U$.

\begin{lemma}
\label{lem:permute}
Let $n = mk+r$, where $0 \leq r \leq k-1$. Let $U \in \gbin{V}{k+r}$ be a negative weight $(k+r)$-dimensional subspace of $V$, and let $\bm{S} \subset \gbin{V}{k}$ be a family of $k$-dimensional subspaces satisfying properties (\lowerromannumeral{1}) and (\lowerromannumeral{2}) of \lref{mrBtotherescue}. If $\pi \in GL(V)$ is an invertible linear transformation that fixes $U$, then
\begin{equation}
\label{spreadpermuted}
\pi(\bm{S}) = \{ \pi(S) : S \in \bm{S} \} \subset \gbin{V}{k}
\end{equation}
is a family of $k$-dimensional subspaces that satisfies properties (\lowerromannumeral{1}) and (\lowerromannumeral{2}) of \lref{mrBtotherescue}. Moreover, some $k$-dimensional subspace in $\pi(\bm{S})$ must have nonnegative weight.
\end{lemma}

\noindent \textbf{Proof.} Since $\pi \in GL(V)$ fixes $U$ and $\bm{S} \subset \gbin{V}{k}$ is a family of $k$-dimensional subspaces that satisfies properties (\lowerromannumeral{1}) and (\lowerromannumeral{2}) of \lref{mrBtotherescue}, we see that $\pi(\bm{S})$ must satisfy the same properties.
Hence,
\begin{equation}
\label{sumSpermute}
\sum_{S \in \bm{S}} b_{\pi(S)} = \sum_{w \in \gbin{V}{1} \setminus \gbin{U}{1}} f(w) = -f(U) > 0,
\end{equation}
so some $k$-dimensional subspace in $\pi(\bm{S})$ must have nonnegative weight. \qed

We show that if $n \geq 3k$ and $T \in \gbin{V}{k}$ is a negative weight $k$-dimensional subspace, then there are almost $\gbin{n-1}{k-1}$ nonnegative $k$-dimensional subspaces that have trivial intersection with $T$.

\begin{lemma}
\label{lem:nonnegT}
If $n \geq 3k$ and $T \in \gbin{V}{k}$ is a negative weight $k$-dimensional subspace, then there are at least
\begin{equation}
\label{nonnegtrivialwithT}
\left(1 - \frac{1}{q^{n-3k+1}} \right) \gbin{n-1}{k-1}
\end{equation}
nonnegative $k$-dimensional subspaces that have trivial intersection with $T$.
\end{lemma}

\noindent \textbf{Proof.} Write $n = mk+r$ where $0 \leq r \leq k-1$. Since $k+r < n$ and $T \in \gbin{V}{k}$ is a negative weight $k$-dimensional subspace, there is a negative weight $(k+r)$-dimensional subspace $U \in \gbin{V}{k+r}$ that contains $T$ by \lref{induct}. Consequently, there exists a family $\bm{S} \subset \gbin{V}{k}$ of $k$-dimensional subspaces satisfying properties (\lowerromannumeral{1}) and (\lowerromannumeral{2}) of \lref{mrBtotherescue}.  Define
\begin{equation}
\label{nonnegdisjointT}
\bm{F} = \left \{ S \in \gbin{V}{k} : b_{S} \geq 0, \; \dim(S \cap U) = 0 \right \}
\end{equation}
to be the family of nonnegative $k$-dimensional subspaces that have trivial intersection with $U$.

Consider a random isomorphism $\pi: V \rightarrow V$ that fixes $U$. For each $k$-dimensional subspace $S \in \bm{S}$, define an indicator random variable $Z_{S}$ by
\begin{equation}
\label{indicator}
Z_{S} := \begin{cases}
         1 & \text{if $\pi(S)$ has nonnegative weight} \\
         0 & \text{otherwise.}
         \end{cases}
\end{equation}
Let $Z = \sum_{S \in \bm{S}} Z_{S}$ and note that $Z \geq 1$ because some $k$-dimensional subspace in $\pi(\bm{S})$ must have nonnegative weight by \lref{permute}. On the other hand, $\mathbb{E}(Z_{S})$ is the probability that a randomly chosen $k$-dimensional subspace that has trivial intersection with $U$ has nonnegative weight. Hence,
\begin{equation}
\label{greatexpections}
\mathbb{E}(Z_{S}) = \frac{|\bm{F}|}{q^{(k+r)k} \gbin{n-k-r}{k}}
\end{equation}
because the denominator counts the number of $k$-dimensional subspaces that have trivial intersection with $U$ by \eqref{useful}. By linearity of expectation,
\begin{equation}
\label{loe}
1 \leq \mathbb{E}(Z) = |\bm{S}| \mathbb{E}(Z_{S}) = \frac{|\bm{S}| |\bm{F}|}{q^{(k+r)k} \gbin{n-k-r}{k}},
\end{equation}
so by \eqref{sizeS}
\begin{equation}
\label{boundF}
|\bm{F}| \geq \frac{q^{(k+r)k} \gbin{n-k-r}{k}}{|\bm{S}|} = \frac{q^{(k+r)k} [k] \gbin{n-k-r}{k}}{q^{k+r} [n-k-r]} = q^{(k+r)(k-1)} \gbin{n-k-r-1}{k-1}.
\end{equation}
By \eqref{disjoint} and \eqref{firstterm} with $a = k+r$ and since $0 \leq r \leq k-1$, we have
\begin{align}
\label{almostall}
|\bm{F}| &\geq q^{(k+r)(k-1)} \gbin{n-k-r-1}{k-1} \geq \left(1 - \frac{1}{q^{n-2k-r}} \right) \gbin{n-1}{k-1} \notag \\
&\geq \left(1 - \frac{1}{q^{n-3k+1}} \right) \gbin{n-1}{k-1}.
\end{align}
Since $U$ contains $T \in \gbin{V}{k}$, each $k$-dimensional subspace in $\bm{F}$ also has trivial intersection with $T$. \qed

Finally, we prove a set analogue of \lref{nonnegT}. The proof of \lref{nonnegTset} is similar to Manickam's and Singhi's proof of \conjref{MMS} when $k|n$ \cite{ms} and has also been observed by others \cite{ahs, frankl, Tyomkyn}.

\begin{lemma}
\label{lem:nonnegTset}
If $\mathcal{T} \in \binom{\cX}{k}$ has negative sum, then there are at least
\begin{equation}
\label{nonnegtrivialwithTset}
\binom{n-2k}{k-1} \geq \left( 1 - \frac{(2k-1)(k-1)}{n-2k+1} \right) \binom{n-1}{k-1}
\end{equation}
nonnegative $k$-element subsets of $\cX$ that are disjoint from $\mathcal{T}$.
\end{lemma}

\noindent \textbf{Proof.} Write $n = mk + r$ where $0 \leq r \leq k-1$. Since $k+r < n$ and $\mathcal{T} \in \binom{\cX}{k}$ has negative sum, an argument similar to that of \lref{induct} shows that there is a $(k+r)$-subset $\mathcal{U} \in \binom{\cX}{k+r}$ with negative sum that contains $\mathcal{T}$. Consider a random permutation $\pi \in S_{\cX}$ that fixes $U$. Partition the $(m-1)k$ elements of $\cX \setminus \mathcal{U}$ into $k$-element sets $\mathcal{S}_{1}, \ldots, \mathcal{S}_{m-1}$, and define the indicator random variable $Z_{i}$ to be $1$ if $\pi(\mathcal{S}_{i})$ has nonnegative sum and $0$ otherwise. Since $0 \leq r \leq k-1$, repeating the first moment method argument of \lref{nonnegT} yields that there are at least
\begin{equation}
\label{firstmomentrepeat}
\binom{n-k-r-1}{k-1} \geq \binom{n-2k}{k-1} \geq \left( 1 - \frac{(2k-1)(k-1)}{n-2k+1} \right) \binom{n-1}{k-1}
\end{equation}
nonnegative $k$-element subsets of $\cX$ that are disjoint from $\mathcal{T}$. \qed

\section{Proof of \tref{vecanalog}}
\label{sec:theproof}

In this section we prove \tref{vecanalog}.

\vspace{0.25cm}

\noindent \textbf{Proof of \tref{vecanalog}} Recall that $A \in \gbin{V}{k}$ denotes the highest weight $k$-dimensional subspace of $V$, and that $C \in \gbin{V}{k}$ denotes the highest weight $k$-dimensional subspace such that $\dim(A \cap C) = 1$. Let $A \cap C = v \in \gbin{V}{1}$. If all $k$-dimensional subspaces containing $v$ have nonnegative weight, then there are at least $\gbin{n-1}{k-1}$ $k$-dimensional subspaces with nonnegative weight in $V$.

Otherwise, some $k$-dimensional subspace $T \in \gbin{V}{k}$ containing $v$ has negative weight. Suppose, for a contradiction, that there are at most $\gbin{n-1}{k-1}$ nonnegative $k$-dimensional subspaces in this case. By \lref{lotson1}, there are at least
\begin{equation}
\label{spacesonv}
\left( 1 - \frac{1}{q^{n-3k}} - \frac{1}{q^{n-2k-1}} - \frac{1}{q^{n-2k}} - \frac{1}{q^{n-2k+1}}  + \frac{q+1}{q^{2n-4k+1}} \right) \gbin{n-1}{k-1}
\end{equation}
nonnegative $k$-dimensional subspaces containing $v$.

Since $T \in \gbin{V}{k}$ has negative weight, by \lref{nonnegT}, there are at least
\begin{equation}
\label{disjointT}
\left(1 - \frac{1}{q^{n-3k+1}} \right) \gbin{n-1}{k-1}
\end{equation}
nonnegative $k$-dimensional subspaces that have trivial intersection with $T$.

Since $T$ contains $v$, none of the $k$-dimensional subspaces counted in \eqref{disjointT} contain $v$. Summing \eqref{spacesonv} and \eqref{disjointT}, there are at least
\begin{align}
\label{grandfinale}
&\left( \left( 1 - \frac{1}{q^{n-3k}} - \frac{1}{q^{n-2k-1}} - \frac{1}{q^{n-2k}} - \frac{1}{q^{n-2k+1}}  + \frac{q+1}{q^{2n-4k+1}} \right) + \left(1 - \frac{1}{q^{n-3k+1}} \right) \right) \gbin{n-1}{k-1} \notag \\
& \geq \left( 2 - \frac{1}{q^{n-3k}} - \frac{1}{q^{n-3k+1}} - \frac{1}{q^{n-2k-1}} - \frac{1}{q^{n-2k}} - \frac{1}{q^{n-2k+1}}  + \frac{q+1}{q^{2n-4k+1}} \right) \gbin{n-1}{k-1}
\end{align}
nonnegative $k$-dimensional subspaces in $V$. For $k \geq 3$ and $n \geq 3k+1$, however, the expression in \eqref{grandfinale} is greater than $\gbin{n-1}{k-1}$, which contradicts our assumption. When $k=2$, we can use \eqref{lotson1k2} in place of \eqref{spacesonv} to see that there are at least
\begin{equation}
\label{grandfinalek2}
\left( 2 - \frac{1}{q^{n-6}} - \frac{1}{q^{n-5}} - \frac{1}{q^{n-3}} + \frac{q+1}{q^{2n-7}} \right)[n-1]
\end{equation}
nonnegative $2$-dimensional subspaces in $V$. For $n \geq 7$, however, the expression in \eqref{grandfinalek2} is greater than $[n-1]$, which contradicts our assumption.

Finally, Manickam and Singhi's \cite[Theorem 3.1]{ms} verification of \conjref{vecanalog} in the case that $k|n$ proves \tref{vecanalog} when $n = 3k$. \qed

\section{Proof of \tref{quadratic}}
\label{sec:thequadraticproof}

In this section we prove \tref{quadratic}.

\vspace{0.25cm}

\noindent \textbf{Proof of \tref{quadratic}} If all $k$-element subsets containing $x_{1}$ have nonnegative sum, then there are at least $\binom{n-1}{k-1}$ $k$-element subsets of $\cX$ with nonnegative sum.

Otherwise, some $k$-element subset $\mathcal{T} \in \binom{\cX}{k}$ containing $x_{1}$ has negative sum. Suppose, for a contradiction, that there are at most $\binom{n-1}{k-1}$ nonnegative $k$-element subsets in this case. By \lref{lotson1sets}, there are at least
\begin{equation}
\label{spacesonx1}
\left(1 - \frac{(6k-3)(k-1)}{n-2k+1} \right) \binom{n-1}{k-1}
\end{equation}
nonnegative $k$-element subsets containing $x_{1}$ since $n \geq 8k^{2}$.

Since $\mathcal{T} \in \binom{\cX}{k}$ has negative sum, by \lref{nonnegTset}, there are at least
\begin{equation}
\label{disjointTset}
\left(1 - \frac{(2k-1)(k-1)}{n-2k+1} \right) \binom{n-1}{k-1}
\end{equation}
nonnegative $k$-element subsets of $\cX$ that have trivial intersection with $\mathcal{T}$.

Since $\mathcal{T}$ contains $x_{1}$, none of the $k$-element subsets counted in \eqref{disjointTset} contain $x_{1}$. Summing \eqref{spacesonx1} and \eqref{disjointTset}, there are at least
\begin{align}
\label{grandfinaleset}
\left( 2 - \frac{(8k-4)(k-1)}{n-2k+1} \right) \binom{n-1}{k-1}
\end{align}
nonnegative $k$-element subsets in $\cX$. For $n \geq 8k^{2}$, however, the expression in \eqref{grandfinaleset} is greater than $\binom{n-1}{k-1}$, which contradicts our assumption. \qed

\medbreak
\noindent{\sc Acknowledgement:}  We are grateful to Simeon Ball, Chris Godsil, Po-Shen Loh, Karen Meagher, Bruce Rothschild, Benny Sudakov, Terence Tao, Jacques Verstra\"{e}te, and Rick Wilson for their advice.

\bibliographystyle{plain}
\bibliography{MMSVecSet}

\end{document}